\documentclass[12pt]{article}
\usepackage{amsmath}
\usepackage{graphicx}
\usepackage{graphics}
\usepackage{algorithm}
\usepackage{algpseudocode}
\usepackage{algorithmicx}
\usepackage{natbib}
\usepackage{url} 
\usepackage{amsfonts,amsthm,amsmath,amsthm}
\usepackage{hyperref,natbib}
\usepackage{setspace}
\usepackage{graphicx}
\usepackage{hyperref}
\usepackage{bm}
\usepackage{autobreak}
\usepackage[normalem]{ulem}
\usepackage[usenames,dvipsnames,svgnames,table]{xcolor}

\newcommand{\X}{\mathcal{X}}
\newcommand{\W}{\mathbf{W}}
\newcommand{\R}{\mathbb{R}}
\newcommand{\Xset}{\mathcal{X}}

\newcommand{\x}{\mathbf{x}}
\newcommand{\y}{\mathbf{y}}
\newcommand{\z}{\mathbf{z}}
\newcommand{\K}{\mathbf{K}}
\newcommand{\A}{\mathbf{A}}
\newcommand{\Y}{\mathcal{Y}}
\newcommand{\nvar}{d}

\newcommand{\veck}{\mathbf{k}}
\newcommand{\vecy}{\mathbf{y}}

\DeclareMathOperator*{\argmin}{argmin}
\DeclareMathOperator*{\argmax}{argmax}

\usepackage{mathtools}

\newcommand{\blind}{0}

\begin{document}
	
	\def\spacingset#1{\renewcommand{\baselinestretch}%
		{#1}\small\normalsize} \spacingset{1}

	
	\if0\blind
	{
		\title{\bf A survey on high-dimensional Gaussian process modeling with application to Bayesian optimization}
		\author{
		Micka\"el Binois\thanks{Universit\'e C\^ote d'Azur, Inria, CNRS, LJAD, France
		Corresponding author: 
				\href{mailto:mickael.binois@inria.fr}{\tt mickael.binois@inria.fr}}
		\and
		Nathan Wycoff\thanks{
				Current institution: Massive Data Institute, McCourt School of Public Policy, Georgetown University, Washington, DC.}
		}
		\maketitle
	} \fi
	
	\if1\blind
	{
		\bigskip
		\bigskip
		\bigskip
		\begin{center}
			{\LARGE\bf A survey on high-dimensional Gaussian process modeling with application to Bayesian optimization}
		\end{center}
		\medskip
	} \fi
	
	\bigskip
	\begin{abstract}

    Bayesian Optimization, the application of Bayesian function approximation to finding optima of expensive functions, has exploded in popularity in recent years.
    In particular, much attention has been paid to improving its efficiency on problems with many parameters to optimize.
	This attention has trickled down to the workhorse of high dimensional BO, high dimensional Gaussian process regression, which is also of independent interest.
	The great flexibility that the Gaussian process prior implies is a boon when modeling complicated, low dimensional surfaces but simply says too little when dimension grows too large.
	A variety of structural model assumptions have been tested to tame high dimensions,
	from variable selection and additive decomposition to low
	dimensional embeddings and beyond. Most of these approaches in turn require modifications of 
	the acquisition function
	optimization strategy as well. Here we review the defining structural model
	assumptions and discuss the benefits and drawbacks of these approaches in
	practice.
		 
	\end{abstract}
	
	\noindent%
	{\it Keywords:} Black-box optimization, low-intrinsic dimensionality, additivity, variable selection, active subspace

	\spacingset{1.5}

\section{Introduction}

A large number of parameters (high dimensionality) is regularly mentioned as a critical challenge for black-box
optimization. The root of the issue is the so-called curse of dimensionality,
as coined by \citet{bellman1966dynamic}: the exponential dependence of complexity on input dimension.
This is all the more true with limited
evaluation budgets, where relying on a surrogate (or metamodel) is a common
practice. 
In BO, the number of variables impacts both the GP surrogate (via the behavior of distances in high dimension, see \citet{Aggarwal2001}), as well as the search for the next designs (as it affects the acquisition function, which is employed in BO in the search for the next design(s) to query).
This dimensionality-scaling difficulty is mentioned in most existing reviews, see for instance
\cite{Shan2010a}, \cite{Viana2014}, \cite{Shahriari2016},
\cite{frazier2018bayesian}, \cite{Ginsbourger2018} or \cite{Stork2020}. This
article is thus incremental, focusing on more recent trends over several
communities (e.g., engineering, operations research, and machine learning) but
is by no means exhaustive.

We focus on Gaussian process (GP) surrogates here, whose popularity originates
from their modeling flexibility and appealing uncertainty quantification (UQ)
properties. Other alternatives could be entertained, but probably not without giving
 up some amount of either flexibility or small sample efficiency.
For examples, see the use of tree models like random forest regression \citep{hutter2011sequential}, tree-structured Parzen estimators \citep{Bergstra2011}, and Bayesian additive regression trees \citep{chipman2012sequential}, spline models with Bayesian adaptive splines \citep{Francom2019}, or neural network surrogate models like Bayesian neural networks \citep{snoek2015scalable} and deep GPs \citep{damianou2013deep} for optimization, e.g., as by \cite{hebbal2019multi}.
Radial basis function interpolation (RBF), closely related to GPs, is also
popular in this setting \citep{Regis2013}, but is comparatively lacking in UQ capabilities.

In BO, the number of variables impacts both the GP surrogate (via the behavior
of distances in high dimension, see \cite{Aggarwal2001}), as well as the search for
the next designs with the acquisition function (by complicating its optimization).
Specifically, the number of
design points required to keep the same quality of approximation grows exponentially
with the number of variables and the volume concentrates on the boundary of
the search space. As a result, designs are comparatively dispersed in high
dimension. We refer to
\cite{sommerville1958introduction,artstein2015asymptotic} for more details on
high dimensional geometry, and to \cite{Zhigljavsky2021} for a discussion more oriented to BO. 
Optimization of the acquisition function to select the next design point benefits from (comparatively) quick evaluation and gradient evaluation, but still suffers from these same high-dimensional effects.
Depending on the difficulty of the
optimization problem, these issues can occur when reaching ten variables, or
manifest for dozens of them. Still, dedicated methods have been shown to work
for billions of variables \citep{Wang2013}-- though under rather limiting assumptions. 

Indeed, stronger structural model assumptions are needed as the dimension increases,
with three main categories. One idea is to reduce the dimension by removing
variables with little or no impact on the output. Another is to define a few
new variables based on linear or non-linear combinations of the original ones.
The last direction is to assume additivity of effects of variables, or groups
of them. With more structure to learn, model inference becomes harder,
adding estimation risk to the difficulties above.
Indisputably, estimation risk can be reduced by avoiding estimation, and this has been
the approach taken by researchers relying on randomly defined structures (though obviously this is at the cost of variance and bias possibly gained from the randomization procedure).
If we want to infer the structure, however, we also must face the fact that in the sequential design framework, only limited data is at first available to estimate characteristics of the function. Either the structure can be learnt and updated ``online" as the data are observed, or we can break the sequential process into two consecutive stages, the first of which is focused on estimating structure and the second on exploiting it. In this latter case, the optimal balance of budget to dedicate to each stage can be quite problem dependent.

The chosen structure affects the acquisition phase in several ways. It can be
exploited for the modeling effort, as in reduced dimension or additivity assumptions, or it can constrain the
optimization domain to aid the acquisition function search. Independently, strategies such as deploying trust
regions as in \cite{eriksson2019scalable} are also available to limit the
effect of the curse of dimensionality by limiting the volume of the search region.
The resulting optimization approach is more local, but is complemented by restarts for globalization.

There exists many additional challenges and refinements for GPs and BO that are out of scope here, for which the corresponding techniques may need to be adapted for high dimension. One such issue is the scaling in terms of number of design points. 
Logically, more designs are needed to learn in larger dimensions, but one can not hope to match the exponential dependence on the dimension. 
Techniques for coping with large data have also attracted a lot of attention, see for instance \cite{hensman2013gaussian,heaton2019case,wang2019exact,kleijnen2020prediction}, where some are independent of the input dimension as is the case for, say, local models. 
In general we will assume here that running the black-box remains limiting compared to running the BO framework.
Other refinements for GP and BO to cope with complex noise modeling (e.g., non Gaussian noise, input dependent variance) or non-stationarity could be adapted, but high-dimension exacerbates the difficulty of the learning task. 
While we focus on unconstrained optimization, batch (or parallel) optimization, constrained optimization and others could be entertained as well. We refer the interested reader to \cite{garnett2022bayesian} for a broader and more introductory overview of BO.

The remainder of this paper is as follows. First, key notions and notations are
introduced in Section \ref{sec:back}. Next, structural assumptions for high
dimensional GP modeling are detailed in Section \ref{sec:highGP}, before consequences
and adaptations for the acquisition function optimization are presented in Section \ref{sec:highBO}. Section \ref{sec:tests} includes a
list of possible test functions. Finally, some practical guidelines and a summary on promising research directions is given in
Section \ref{sec:conc}.

\section{Background}
\label{sec:back}

Let us consider an expensive-to-evaluate \emph{black-box} simulator $f : \X \subset \R^d \rightarrow \R$ that we want to globally optimize:
\begin{equation}
\text{find~} \x^{*} \in \argmin \limits_{\x \in \X} f(\x).
\label{eq:opt_pb}
\end{equation}
By black-box, we mean that nothing is assumed known about the functional form of $f$: $f(\x)$ can only be queried at any given input point $\x$ (sometimes gradients are also assumed available).

\subsection{The Gaussian stochastic process}

Given an index set $\mathcal{X}$, in our case typically a closed and bounded subset of $\mathbb{R}^\nvar$, a \textit{stochastic process} is simply a rule for assigning to $B$ members of that set $\x^{(1)}, \x^{(2)}, \ldots, \x^{(B)}$ a joint distribution of random variables; maybe it has a density $\delta \left(y(\x^{(1)}), y(\x^{(2)}), \ldots, y(\x^{(B)}) \right)$ with respect to some dominating measure. If the joint distribution is assumed to be Gaussian, we simply need a way of deciding what the mean vector and covariance matrix are going to be for any possible set of $B$ points. These are naturally referred to as the \textit{prior mean function} $\mu$ and \textit{covariance function} $k$ (also \textit{kernel function}). Such a mathematical construct taken all together is called a \textit{Gaussian process} (GP), see, e.g., \cite{Rasmussen2006}.

If we observe $y(\x^*)$ (possibly corrupted by Gaussian noise), we can imagine that there is some stochastic process that maps $\x ^{(1)}, \x^{(2)}, \ldots, \x^{(B)}$ to $\delta \left(y(\x^{(1)}), y(\x^{(2)}), \ldots, y(\x^{(B)}) | y(\x^*) \right)$, that is, that maps any set of points to its conditional distribution given $y(\x^*)$. Happily, this new stochastic process also happens to be a GP (conjugacy with respect to the normal likelihood), and its new mean and covariance functions are available in closed form. In particular, if we observed $y(\x^{(1)}), \ldots, y(\x^{(n)})$, then:

\begin{align*}
m_n(\x) &= \mu + \veck(\x)^\top \K^{-1} (\y - \mu \mathbf{1}),\\
s_n^{2}(\x) &= k(\x, \x) - \veck(\x)^\top \K^{-1} \veck(\x)
\end{align*}
where $\vecy:= (y_1, \dots,
y_n)$, $\veck(\x) := (k(\x, \x^{(i)}))_{1 \leq i \leq n}$, $\K:= (k(\x^{(i)},
\x^{(j)}) + \tau^2 \mathbf{1}_{i=j})_{1 \leq i,j \leq n}$. $\tau^2$ is the noise
hyperparameter, when assuming $y_i = f(\x^{(i)}) + \varepsilon_i$, with
$\varepsilon_i = \mathcal{N}(0, \tau^2)$.

Often, the prior mean function $\mu$ is chosen to be zero (after centering observed $y$): the GP's focus is on defining and exploiting spatial covariance. However, it can definitely be useful to have different $\mu$. If the objective function is known to be nonstationary, adding in polynomial trend terms can be a way to tame this. However, if the proper trend is not \textit{a priori} known, it pays to be cautious when using one in high dimension as it's easy to accidentally extrapolate rather severely and perniciously \citep{Journel1974}. Sparse trend selection might be preferable, as
performed for instance by \cite{kersaudy2015new} in a combination with
polynomial chaos expansion, as may be using basis elements informed via cross-validation
\citep{Liang2014} or a Bayesian framework \citep{Joseph2008}.

The choice of covariance function family determines the qualitative properties of the GP. Of course, we want to make sure that when we're done applying the kernel function to all our data, the covariance matrix that comes out is positive semi-definite. Functions which guarantee this are themselves called positive semi-definite. 
Attention is often further restricted to \textit{stationary} kernels, where the kernel function $k$ is a function of $\x - \x'$, $k(\x,\x')=\tilde{k}(\x-\x')$. 
Within this class, the most popular kernel is probably the squared exponential one in product form, $k(\x, \x') = \sigma^2  \prod \limits_{i = 1}^d \exp \left(-\frac{(x_i-x'_i)^2}{2\theta_i} \right)$. This kernel function is defined only up to its free parameters $\theta_1, \dots, \theta_d$, the lengthscales which determines at what distance covariance begins to drop off between points, and $\sigma^2$ which converts correlation to covariance based on the scale of $y$.  This kernel induces an infinitely differentiable metamodel, which may or may not be desirable, depending on how well that matches our understanding of the black-box. For precise control over differentiability\footnote{For more details on what exactly is meant by the derivative of a stochastic process see \cite[Chapter 10]{Papoulis1965random}.}, consider the Mat\'ern class indexed by smoothness parameter $\nu$, which is defined in terms of the modified Bessel function of the second kind (though can be written in terms of simple functions when $\nu+0.5$ is a whole number). See \citet[Chapter 4]{Rasmussen2006} for more on the Mat\'ern and other covariance function classes. Such GPs can be differentiated at most $\nu-0.5$ times. 

This dependence of GPs on distances between points will be problematic in high dimension. And, consequently, many of the approaches we review in this article will define a kernel function not on the original input space, but on some transformation thereof. Of course, it is possible to view the transformation and the kernel function thereon together as a novel kernel function. But we will instead take the perspective in this article of separating dimension-reducing transformations and the subsequent purely distance-based kernel, as we view this modularity as helpful in model design and comparison.

Insofar as a Gaussian process is the infinite dimensional generalization of the multivariate Gaussian distribution, it retains many of the analytical properties that have propelled Gaussianity to its prominence. Of particular use to us will be the closed form evaluation of certain exceedance probabilities and related metrics which form the basis of deploying these Gaussian processes to optimization problems.

We refer to \cite{Rasmussen2006,Gramacy2020} for
additional details on GP regression, and, e.g., to \cite{kanagawa2018gaussian} for a
discussion of connections to other kernel methods.

\subsection{Bayesian optimization}

Popularized with the Efficient Global Optimization method of \cite{Jones1998},
Bayesian optimization, see, e.g., \cite{Mockus1978}, relies on the GP
probabilistic prior of $f$ based on a few initial observations to define an acquisition function $\alpha: \Xset
\to \R$. This is later optimized to select new points to evaluate
sequentially. 
The initial stage, known as the Design of Experiments (DoE), typically starts by evaluating designs based on optimized Latin hypercube samples, for their space filling properties, see, e.g., \cite{Gramacy2020}.
A pseudo-code is given in Algorithm \ref{alg:bBO}.

\begin{algorithm}[!ht]
\caption{Pseudo-code for BO}\label{alg:bBO}
\begin{algorithmic}[1]
\Require  $N_{max}$ (total budget), GP model trained on initial DoE of size $n = n_0: (\x^{(i)}, y_i)_{1 \leq i \leq n}$
\While{$n \leq N_{max}$}
    \State Choose $\x^{(n+1)} \in \arg \max_{\x \in \X} \alpha(\x)$
    \State Update the GP model by conditioning on $(\x^{(n+1)}, y_{n+1})$
    \State $n \leftarrow n + 1$
\EndWhile
\end{algorithmic}
\end{algorithm}

The expected improvement \citep[EI]{Mockus1978} and upper confidence bound \citep[UCB]{srinivas2010gaussian} criteria are typically used
due to their analytical tractability in both value and gradient. Other choices are possible, see e.g., \cite{Shahriari2016}.
The generic acquisition function optimization problem is:
\begin{equation}
\text{find~} \x^{*} \in \argmax \limits_{\x \in \X} \alpha(\x).
\label{eq:opt_acq}
\end{equation}

Some acquisition functions are natively suitable for when the black-box is
noisy, while some need adaptation like EI, see e.g., \cite{Letham2018} for a
discussion. A pinch of noise might be beneficial for its regularization
effect, see, e.g., \cite{Gramacy2012a}, but the effort of finding a solution
increases drastically as the signal to noise ratio sinks. Not to mention
input-varying and non-Gaussian noises that are hard to deal with even in low
dimension. We refer to \cite{Forrester2008,Roustant2012,Gramacy2020,garnett2022bayesian} for more
details on how to apply BO in practice.

\subsection{The curse of dimensionality and its effects}

\begin{figure}[htpb]
    \centering
    \includegraphics[scale=1.0,trim={1cm 1cm 0.2cm 0.2cm},clip]{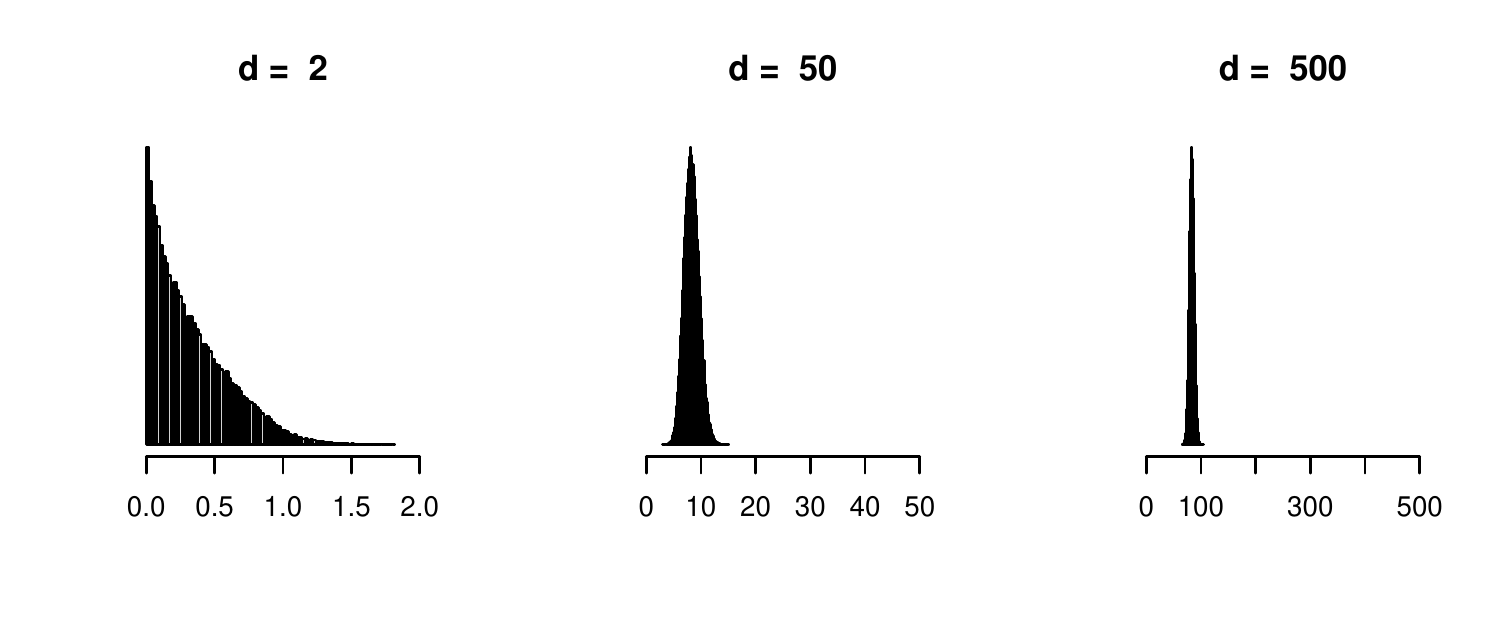}
    \caption{\textbf{Distances Concentrate in High Dimension:} Randomly sampling 10,000 points according to the uniform measure in $[0,1]^d$ and then calculating squared inter-point distances reveals that these concentrate within the bounds of possible values in high dimension: $[0, d]$.
    }
    \label{fig:dists}
\end{figure}

An increase in input dimension is first felt in distance computation, which is at the core of most common covariance kernels.
The difficulty is that points are relatively further away in high dimension, such
that it becomes harder to learn with a covariance based on distances. As is illustrated in Fig.\ \ref{fig:dists}, for large $d$, the squared distance between uniformly sampled points concentrates increasingly further away. In particular, the expected squared distance approaches to $\frac{\nvar}{6}$ with standard deviation relative to the width of the interval of possible distances $[0,\sqrt{d}]$ decreasing as $\sqrt{\frac{7}{180\nvar}}$ \footnote{The results follow from the fact that the density of the squared distance $d$ between two uniform variables on $[0,1]$ is given by $\frac{1}{\sqrt{d}}-1$.}. 
See \cite{koppen2000curse} for a discussion of other implications and properties of high dimensional cubes.

Except for isotropic kernels, the second effect is the increase of the number
of hyperparameters of the models, typically with $d$ many lengthscale parameters to
tune for product covariance kernels. Maximizing the (non-convex) likelihood thus becomes
increasingly hard, since the curse of dimensionality also affects the gradient
based optimization routines typically used. Grouping variables to share the same
lengthscales alleviate this, which has been done either manually based on problem-specific
considerations as done by \cite{Binois2015c} or using the Bayesian information
criterion (BIC) as done by \cite{blanchet2017specific}.

Similarly, accurately optimizing the acquisition function is complex in high
dimension, especially since it is a multi-modal optimization task, with flat
regions. Furthermore, given that most of the volume is on the boundary of the
domain and that GP-based predictive variance increases with the distance to
the designs -- a desirable property at least in low dimension -- the downside
is that the optimum is generally found on one of the exponentially many
vertices or sides of the $d$-hypercube, tilting the balance towards blind
exploration.

The curse of dimensionality also manifests itself when random sampling,
required for Monte Carlo estimation of certain integral quantities, among them entropy, the optimizing point of a posterior realization (as is required in Thompson sampling), and certain global sensitivity measures. The issue is that a uniformly sampled point will be increasingly distant, on average, to any given point as dimension increases.
Hence the challenge is not only to
scale in terms of modeling accuracy but also to keep inference manageable and
avoid pitfalls in acquisition function optimization.

\section{High-dimensional Gaussian process modeling}
\label{sec:highGP}

Structural model assumptions are the only way to avoid the exponential dependence on
the dimension in modeling, for which various options are summarized for
instance by \citet{bach2017breaking}, who also gives the corresponding
generalization bounds.
These different structural model assumptions have been
adapted in the GP context, plus others, as we detail next. One exception is to
use isotropic kernels and thus always a single lengthscale, for which increasing the dimension influences only the behavior of distances.

\subsection{Variable selection or screening}

Whenever possible, limiting the number of variables based on expert knowledge
on the optimization problem is recommended. If expert knowledge is unavailable,
a simple, data-driven idea is to
perform variable selection, or screening, before optimization, e.g.,
using the Morris technique \citep{morris1991factorial} or via hierarchical
diagonal sampling as in \cite{Chen2012}. Other global sensitivity analysis
techniques can be used to select variables individually and we refer to
\cite{Iooss2015} for an entry point to this topic.

Hence, one of the early attempts to tackle high-dimensions is to assume that
most of the variables have no effect: 
\begin{equation}
    \text{model:~} f(\x) \approx g(\x_I) \text{~with~} I \subset \{1, \dots, d \}, |I| \ll d
    \label{eq:varsel}
\end{equation}
and then identify those
influential variables in the set $I$. 
For the Gaussian and other
stationary kernels in tensor product form, this can be performed by looking at
the lengthscale values: the
covariance varies less for very large values of the lengthscale (when parameterized as above), whose corresponding input variables can thus be removed. 
This is also known under the term automatic relevance
determination (ARD), as discussed, e.g., by \citet[Chapter 5]{Rasmussen2006}.
\cite{Salem2018} shows that this indeed holds asymptotically.
Hence the idea is to rank the variables based on their lengthscales to determine
the more influential from the less. These influential variables are then
used to build a GP model and optimize expected improvement.

Nevertheless, $\theta_1 > \theta_2$ does not necessarily imply that $x_1$ is
less important than $x_2$, see for example \cite{wycoff2021sequential,lin2020transformation}.
Hence rather than looking only only at the lengthscales values, \cite{linkletter2006variable} propose to compare the posterior distributions of lengthscales corresponding to both real and artificially added inert variables in a fully Bayesian GP framework. 
This is complemented by a local variable selection in \cite{Winkel2021}, relying of predictions made when ignoring some dimensions. 
\cite{Eriksson2021} choose instead to put a Horseshoe prior \citep{Carvalho2009} on the inverse lengthscales and perform gradient-based numerical posterior simulation in order to sample from this posterior of functions defined on subsets of variables.
With a more incremental approach, \citet{Marrel2008} propose to sequentially add
variables to the regression (resp.\ covariance) element based on the Akaike
information criterion (resp.\ predictivity coefficient $Q_2$).

Rather than definitively deciding on which variables to keep, which can be
time-consuming, \citet{ulmasov2016bayesian} propose to sample a few variables
at each iteration, where the weight vector is determined with principal
component analysis (PCA) on the $(\x)_{1 \leq i \leq n}$.
This may not be efficient for small budget as
indicated in \citet{Li2017}, who prefer to select variables uniformly and to
fill in values for the remaining variables.

In practice, most variables typically have a possibly limited but non-null influence on the output.
Though this can mean that choosing the number of variables to keep is somewhat arbitrary,
the simplicity, both in terms of implementation and interpretation, of variable selection makes it a 
compelling approach to dimension reduction, particularly for a first pass.
Of course, the variable selection itself becomes only harder as the dimension increases, but this is
to a lesser extent than some of the more sophisticated approaches discussed below.
Nevertheless, the approach's strength is also its limitation, and it is common for the complicated functions which 
make up black-box problems to truly and strongly depend on all input parameters.
In the case where all variables have the
same influence, variable selection would even fail to reduce the dimension by one no matter how simple the relationship.
Other structural assumptions can overcome these limitations, as well as exploit interactions between variables.

\subsection{Additive and ANOVA models}

One set of structural assumptions amenable to keeping all variables but limiting their interaction is that of additivity: 
\begin{equation}
\text{model:~} f(\x) \approx \mu + \sum \limits_{i = 1}^d g_i(x_i)
\label{mod:add}
\end{equation}
with univariate functions $g_i$.
This has been transposed to the GP framework by
\cite{neal1997monte,plate1999accuracy,Durrande2010,Duvenaud2011}, initially via
the summation of univariate kernels: $k(\x, \x') = \sum \limits_{i = 1}^d
\sigma_i^2 k_i(x_i, x'_i)$, which results in a valid covariance function just as the product form. A useful property for interpretability and visualization is
that the GP predictive mean can be decomposed into a sum of univariate
components: $m_n(\x) =
\veck(\x)^\top
\K^{-1} \vecy = \sum \limits_{i = 1}^d \veck_i(x_i) \K^{-1} \vecy = \sum \limits_{i = 1}^d m_{n,i} (x_i)$ with
$\veck_i(x_i):= (k_i(x_i, x^{(j)}_i))_{1 \leq j \leq n}$. A more surprising
property is that the covariance can become non invertible due to linear
relationships appearing between distinct observations -- see e.g., \cite{Durrande2010}
for examples--which can be alleviated by adding a noise term. Consequently,
the predictive variance can be zero at an unobserved design point, a somewhat
detrimental side effect for exploration and thus optimization. Another difficulty is the estimation
of the hyperparameters, with the need of an additional variance parameter $\sigma_i$ per
coordinate (thus approximately doubling the number of hyperparameters). Still, compared to the tensor product form whose values quickly
go to zero in high dimension, the sum form scales much better. Making the
black-box more additive by applying a transformation of the outputs is explored by \cite{lin2020transformation}.

Higher-order models can be defined in the same way \citep{Duvenaud2011},
usually restricted to order two or selecting higher order components only if
all lower orders are already selected. 
\cite{plate1999accuracy} instead proposes adding all interactions at once on top of first order ones.
Directly identifying groups of variables is also possible, see e.g.,
\cite{Kandasamy2015,Gardner2017,Wang2017,Wang2018}:
\begin{equation}
\text{model:~} f(\x) \approx \mu + \sum
\limits_{i = 1}^M g_i(\x_{A_i}) 
\label{mod:add_block}
\end{equation} 
with multivariate $g_i$ functions acting on $A_i$ disjoint subsets of $\{1,
\dots, d\}$ such that $\bigcup_i^M A_i = \{ 1,
\dots, d \}$. The restriction of disjoint subsets of variables have been
further lifted in subsequent works, see e.g.,
\cite{rolland2018high,hoang2018decentralized}.

Similar in their form but rooted in global sensitivity analysis are additive
models based on the functional analysis of variance (fANOVA, a.k.a.\ Sobol-Hoeffding) decomposition
\citep{efron1981jackknife,sobol2001global}: 
\begin{equation}
\text{model:~}f(\x) \approx c +
\sum \limits_{i = 1}^d g_i(x_i) + \sum \limits_{j < k} g_{jk}(x_j , x_k) +
\dots + g_{12 \dots d}(x_1, x_2, \dots, x_d)
\label{mod:anova}
\end{equation} with elementary functions $g_{\dots}$ that are required to be centered and orthogonal for the uniqueness of the decomposition.
The benefit is that sensitivity analysis can be conducted and
interpreted as in a regular ANOVA. \cite{Muehlenstaedt2012} rely on this
decomposition up to second order interaction to build their model,
choosing their components on the basis of Sobol indices that appear in this formulation
(and estimated by a first-pass GP with an anisotropic tensor-product kernel).
These are the so-called main effects (one variable only) and total
interaction effects (effect of two variables at any order). Selecting which
interaction to remove to form cliques requires a thresholding scheme. 
\cite{ulaganathan2016high} proposes a similar approach with the addition of cut-points when gradient observations are available. 
\cite{Durrande2013} go further with a dedicated kernel whose form is
$k_{\mathrm{ANOVA}}(\x, \x') = \prod \limits_{i=1}^D (1 + k^i(x_i, x'_i))$, as
introduced by \cite{stitson1999support}. There the sensitivity indices of the
ANOVA representation are analytically tractable. In \cite{Ginsbourger2016}, the
ANOVA decomposition is performed directly on the kernel and propagated to the
corresponding random field under appropriate orthogonality conditions. Looking for
sparsity to avoid estimating the entirety of the $2^d$ components, they define
projectors that allow separation of the additive components (with
cross-covariance between main effects) from their complement.

The main downside of these techniques is that inference is challenging in this
context, with combinatorially many terms. Hence various techniques to estimate
the hyperparameters have been applied: coordinate-ascent-like
\citep{Durrande2010} or quasi-Newton \citep{Duvenaud2011} methods for likelihood
maximization. Relying on randomness is sometimes preferred to bypass the cost
of a full optimization, like by \cite{Kandasamy2015} where random
decompositions are sampled and the best one for the likelihood is selected (fixing
the order and number of terms), or by \cite{Wang2018}. \cite{Gardner2017}
attempts to elicit the structure via a dedicated Metropolis-Hastings
algorithm. For overlapping subsets, \cite{hoang2018decentralized} also rely on
random groups while \cite{rolland2018high} use a dependence graph and Gibbs
sampling to perform inference. 
In this vein of methods, there is \cite{delbridge2019randomly} that
reconstruct the high dimensional kernel by a sum of univariate ones in random
directions, in the spirit of the turning band method from geostatistics
\citep{Journel1974}.

In terms of advantages, on the other hand, this method maintains the 
interpretability of variable selection, particularly if the selected model includes mostly first or
second order terms and is somewhat sparse. Learning which variables interact by observing
which pairs are selected by an inference-driven procedure can be scientifically interesting in and of itself.

Some works are dedicated to scaling to many observations as well, such as
\cite{mutny2018efficient} with a basis expansion of GP kernels.
\cite{Sung2019} also propose a basis expansion, further complementing it with
a multi-resolution scheme with a group lasso estimation procedure.
\cite{Wang2018} scale to many observations with a random partition of the
input space, random additive approximation and random features decomposition
of the kernel.

Here, one underlying hypothesis is that high order interaction components are
negligible (because they are hard to estimate in this formalism), though this is a long-standing tradition in experimental design. 
Another concern is that ANOVA cannot detect non-linearity and multi-modality
\citep{Palar2017}. These downsides are mitigated by the following framework.

\subsection{Linear embeddings}
\label{sec:GP_AS}

One approach to designing high dimensional kernels is to avoid doing so directly,
and instead run the data through a dimension reducing map first. 
The image of these data under this map is called its embedding, and an obvious
class of functions to use for dimension reduction is the class of linear functions.
Indeed, as noted by \cite{Marcy2018}, this idea of using a linear mapping dates back to at least
\cite{matern1960spatial}.
The question is then simply which linear function to use in particular.

Denote this generic mapping by $\z = \A \x$ with $\A \in
\R^{r \times d}$. In the case of \cite{vivarelli1999discovering}, $r=d$, so the mapping
serves to rotate the space rather than actually reduce the dimension, but we are more interested in the $r \ll d$ case.
When $r = 1$, this is a popular dimension reduction technique called the single index
model: 
\begin{equation}
    \text{model:~} f(\x) \approx g(\mathbf{a}^\top \x) \text{~with~} \mathbf{a} \in \R^d.
    \label{mod:single}
\end{equation} 
We refer to \cite{Gramacy2012} for the GP treatment, with $\mathbf{a}$
treated as an additional kernel hyperparameter. $\mathbf{a}$ could instead simply be randomly
chosen at each iteration, as suggested by \cite{kirschner2019adaptive}.
Unfortunately, extending to $r>1$ is by no means trivial. The assumption is:
\begin{equation}
    \text{model:~}f(\x) \approx g(\A^\top \x)
    \label{mod:ridge}
\end{equation}
 where such functions are called
\emph{ridge functions} in the literature. It corresponds to the observation, sometimes
backed by theoretical evidence \citep{Constantine2016}, that the variation of
high dimensional functions can be concentrated around a few but unknown
directions. There are several approaches to choosing $\A$, which we'll review in turn.

For GP regression, arguably the most direct way is simply to treat $\A$ as just another
kernel hyperparameter to learn, e.g., by marginal likelihood optimization. 
Fixing $r$, \cite{Garnett2013} provide an approximately Bayesian
scheme to do so, with a Laplace approximation to the likelihood followed by an approximate
marginalization over hyperparameters. \cite{Tripathy2016} rely on a two stage
approach: first using the likelihood to learn an orthogonal $\A$, then finding
the rest of the hyperparameters (and repeat). BIC is used to identify $r$. The
orthogonality constraint requires some special consideration, since these matrices lie
on the Stiefel manifold. \cite{seshadri2019dimension} reinterpreted 
a similar problem, drawing connections between regression and approximation.
Also building on \cite{Tripathy2016}, \cite{Yenicelik2020} observed that
the likelihood may not always select the best matrix among alternatives. A
full Bayesian treatment was proposed by \cite{Marcy2018}, with priors on
matrix manifolds and the dimension, requiring advanced Monte Carlo methods.

Running a sensitivity analysis independent of the GP model to select $\A$ as a kind
of preprocessing step is also an option.
When gradients of the black-box are available, recovering the matrix $\A$ is relatively
straightforward. It corresponds, up to a rotation, to the eigenvectors with
non-zero eigenvalues of the matrix $\mathbf{C} := \int_\Xset
\nabla(f(\x))^\top
\nabla(f(\x)) \lambda(d\x)$ where $\lambda$ is any well-behaved measure on the design space (typically
the Lebesgue one on hypercubic domains). Then a Monte Carlo estimator can be
used. Loosening the assumption that some eigenvalues will be exactly zero leads to
the active subspace (AS) methodology \citep{Constantine2015}, analyzed for
dimension reduction and visualization. In particular, the presence of an AS is
hinted at by gaps in the eigenvalues of $\mathbf{C}$ \citep{Constantine2015}. If
the gradient can be evaluated, this is a good pre-processing step before GP modeling
\citep{Eriksson2018}. Without the gradient information, finite differencing is
generally too costly, such that using a GP to estimate $\mathbf{C}$ might be
more appropriate, e.g., as in \cite{Fukumizu2014,Palar2017}.
\cite{Djolonga2013} uses directional derivatives with finite difference to
recover $\A$ with low rank matrix recovery before BO. Other compressed sensing
techniques can be used, as done initially by \cite{carpentier2012bandit} and later by \cite{groves2018efficient}.
The MC estimators of $\A$ are applied in a two stages approach, due to the assumption of iid sampling of the design points.
\cite{wycoff2021sequential} showed that an estimate of the $\mathbf{C}$ matrix
of a GP is tractable for standard stationary kernels, relieving the need for this sampling assumption.
Other works in the vein of sensitivity analysis on slice inverse regression or partial least
squares to recover $\A$ and possibly update it during optimization
\citep{Bouhlel2016,Zhang2019,chen2020semi}.  
\cite{Lee2019} advocate for a
modified AS matrix which exaggerates the influence of the average gradient (over the input space).

And finally, rather than learning $\A$ as part of inference on GP hyperparameters
or as the result of a sensitivity analysis, a third option is to simply select $\A$ randomly, 
either by generating a single $\A$ matrix before the modeling step or using different $\A$s, say,
at each iteration.
For example, in random embedding BO (REMBO), \citet{Wang2013,Wang2016} use a fixed and randomly sampled $\A$.
One justification is the
stability of random projection of the $L_2$ norm due to the
Johnson-Lindenstrauss lemma \citep{johnson1984extensions}, as noticed, e.g., by
\cite{letham2020re}. In the context of BO, the rationale is further that, at least
for an unbounded domain, there exists a solution to the problem on the low
dimensional embedding, explaining the success of random optimization
\citep{Bergstra2012} on some hyperparameter tuning problems. 

Even with a candidate $\A$ in hand, whether through optimization, sensitivity analysis, or random sampling, 
we still have modeling decisions to make (not to mention acquisition decisions; Section \ref{sec:BO_AS}).
The seemingly innocuous choice of a hypercubic domain can give headaches when combined with a linear
dimension reduction, as the possible solution space is no longer simply a cube but a polytope defined by $\A$.
If our GP is fit to the low dimensional space, it may have an easier time finding a next
candidate point to optimize, but we have to figure out which point in the original space corresponds to 
this low-dimensional optimum. This means that not only are preimages non-unique, but the backprojection of a point is not even
guaranteed to live in the original unit hypercube, and a convex projection is possibly required
to regain feasibility. This means that our backprojection is no longer linear. 
The remediation options proposed in
\cite{Wang2013} are to ignore the corresponding non-injectivity issue, or to
use a kernel defined on $\X$, losing the benefits of low dimensional GP
modeling. \cite{Binois2015} proposed to include high-dimensional information
in $\Y$ with a warping to address the non-injectivity issue, while
\cite{binois2020choice} defines an alternative mapping from $\Y$ to $\X$ to
avoid it. \cite{nayebi2019framework} also bypasses this problem by choosing a
sparse random matrix with $\left\{-1,0,1
\right\}$ elements only, in essence selecting from embeddings on diagonals of the hypercube.
These ideas can be extended to affine embeddings, as proposed for global
optimization by \cite{cartis2020constrained,cartisdimensionality}.

Then comes the question of how to use high
dimensional data that are not on the linear embedding (unlike for REMBO).
Unless $\A$ is perfectly recovered, introducing a noise term is necessary to account
for the discrepancy. To this end, \cite{moriconi2020high} uses GP regression
on quantiles for axis-aligned projections. Another question is how to select
the lengthscales in this context since the product kernels are not preserved
in the embedded space, for which \cite{letham2020re} show that a specific
parameterization (Mahalanobis distance based kernel) is preferable and handles
distortions. Still, inference is complicated whereas the REMBO program needs
only to fit a low dimensional GP. What distribution to sample $\A$ from has also been
briefly investigated in \cite{Binois2015c,binois2020choice,letham2020re}.

Though in general any GP can be fit in the reduced space, links to additive assumptions are natural,
and \cite{gilboa2013scaling} combine these in the projection-pursuit style, leading to projected additive approximations. 
In this context, tuning the $r$ parameter is performed by adding dimensions iteratively and halting when
the accuracy stops increasing. \cite{Li2016} further discuss how to extend this approach
to optimization, in a restricted-projection pursuit setup. 
This can be viewed as somewhere in between fitting a full low dimensional model and using first order additive model.  
These ideas can be combined by estimating several low-dimensional subspaces,  see e.g., \cite{Yenicelik2020,wong2020embedded}.

The estimation risk is strong when using a linear embedding (which has $p\times r$ parameters in general), with unresolved questions: i) How well
must the matrix $\A$ be learned before it is better than directly fitting the high dimensional problem? ii) How much of the budget should be
dedicated to this task? iii) When is it better to learn $\A$ dynamically, leading
to a noisy GP (whenever design points are not all on the same embedding as for REMBO)? 
Plus finding an appropriate $r$ remains difficult unless expert
information is given. Fortunately, taking larger values than necessary is not
detrimental --except that it quickly reintroduces high-dimensional challenges \citep{Wang2013,cartisdimensionality}. 

Already, the linear embeddings approach jettison much of the interpretability that selection and ANOVA approaches offer (it will get worse when we consider nonlinear embeddings in the next section). 
In the case when $r\in\{1,2\}$, it is possible to create visualizations of the function, which can be great sources of intuition. 
But in larger dimensional cases, we are reduced to squinting at the loadings of each variable in the retained directions.
On the other hand, there are indeed functions which can be reduced to much smaller dimension when using linear embeddings compared to the axis-aligned variable selection \footnote{A toy example is $\sin(\frac{4\pi}{d}\mathbf{1}^\top\mathbf{x})$.}.
Which of these two approaches will be more appropriate is quite problem-dependent.

\subsection{Non-linear embeddings and structured spaces}

Dropping the linearity assumption adds even more flexibility to
the model at the cost of requiring even more data to fit it. Recovering such a manifold
suitable for regression has been proposed e.g., by \cite{Guhaniyogi2016}, when
the data is on a low dimensional manifold.
A simple way to extend linear dimension reduction to the nonlinear case is to use a locally linear approach, as does \citet{wycoff2021gradient} in the context of AS.
Also defined by analogy to AS is \citet{Bridges2019}'s 1D active manifold, not yet applied to GPs.
In this line of thought, it is also possible to include
generative topographic mapping \citep{Viswanath2011}, GP latent variable
models, e.g., \citep{Lawrence2005,titsias2010bayesian}, or deep GPs
\citep{damianou2013deep,hebbal2019multi,sauer2020active}.
Deep Gaussian processes refer to the modeling strategy of assuming that the input is linked to 
the output via a chain of Gaussian processes, the output of one serving as 
the input for the next. This is not to be confused with the concept of deep \textit{kernels} \citep{Wilson2016,huang2015scalable}, 
which involves only a single Gaussian process, but whose kernel function is parameterized by a neural network.
Typically, inference is conducted on the GP and neural net weights via optimization of (an estimate of) the log marginal likelihood, but in the event of structured inputs or partial labeling the neural network can be initialized as an autoencoder (i.e., initialized to reconstruct the input).
An orthogonal direction
to avoid the larger optimization budgets is via multi-fidelity, when cheaper
but less accurate version(s) of the black-box are available, as exploited
e.g., by \cite{Ginsbourger2012,Falkner2018}.

A setup where these highly flexible models may be more amenable is in the case
where some additional information is available about the structure of the
problem. One such case is the optimization of geometric shapes, typically
airfoils. AS is performed in this context by \cite{lukaczyk2014active}. As
many options to parameterize these shapes are available, some are more adapted
to optimization. Independently of this choice, \cite{gaudrie2020modeling}
propose to work in the shape space, defined with shape eigenvectors and their
values, which is not costly since the computation of the geometry is in
general negligible next to the subsequent simulation. A non-linear computation
of the shape basis is used in \cite{chen2020airfoil}, with generative
adversarial networks used to learn a manifold from real data. 

A related approach uses variational autoencoders \citep{Kingma2013} as a deep kernel \citep{Gomez2018}, 
which ensures that the latent state upon which the Gaussian process applies its kernel can
approximately recover the original, unencoded input. 
This approach is popular when optimizing over structured non-Euclidean inputs, such as molecules \citep{Gomez2018,Tripp2020,Maus2022,Deshwal2021,Grosnit2021,Eissmann2018}, for which standard BO is unavailable.
Additional examples on functional data indexed on graphs include \cite{espinasse2014parametric}.

\cite{jaquier2020high} perform geometry-aware BO on Riemannian manifolds for
robotics, transposing some of the ideas from Section \ref{sec:GP_AS} to
non-Euclidean spaces.
A somewhat related assumption by \cite{oh2018bock} is to use cylindrical
coordinates instead of the original Cartesian ones (and rather than polar
coordinates as \cite{padonou2016polar}, that do not scale well). The underlying
assumption is that the solution is close to the center of the domain, if
suitably chosen. Transforming the coordinates amounts to separating radial and
angular components: $$T(\x) = \left\{
\begin{array}{ll}
(||\x||_2, \x / ||\x||_2) & \text{for~} ||\x||_2 \neq 0 \\
(0, \mathbf{a}_{arbitrary}) & \text{for~} ||\x||_2 = 0
\end{array}
\right.,
$$ while the inverse transformation is $T^{-1}(r, \mathbf{a}) = r \mathbf{a}$.
The corresponding covariance kernel is $k_{cyl}(\x, \x') = k_{r}(r, r') \times
k_{\mathbf{a}}(\mathbf{a}, \mathbf{a}')$. The 1D radius component $k_{r}(r,
r') = k( |(1-r^\alpha)^\beta - (1-r'^\alpha)^\beta| )$, $\alpha, \beta > 0$ is
chosen to further focus on the center \citep{oh2018bock} while the angular
component is a continuous radial kernel $K_{\mathbf{a}}(\mathbf{a},
\mathbf{a}') = \sum \limits_{p = 0}^P c_p(\mathbf{a}^\top \mathbf{a}')^p$,
$c_p > 0$, $\forall p$ \citep{jayasumana2014optimizing} with a user defined
$P$.

The nonlinear approaches have seen great success in high dimension, particularly when the space is structured/non-Euclidean.
However, this comes at the cost of added complexity, both in computation (e.g., when using deep neural networks) and interpretation.
Therefore, our recommendation is to reserve these techniques for the toughest problems, and consider
simpler solutions to simpler ones, particularly if insight into function, and not only an optimal solution, is desired.

There are many hoops and hurdles to construct a satisfying GP model in
high dimension. Yet, this only gets one half the way, since it remains to
optimize the acquisition function.

\section{High-dimensional acquisition function optimization}
\label{sec:highBO}

A proposed high-dimensional GP model is usually combined with a strategy for
acquisition function optimization, adapted to the specific model structure and the curse of dimensionality.
Many of these approaches can be extricated from the modeling framework in which they were proposed, leaving many possible combinations unexplored. We now discuss them on their own merits.

Faster evaluation time and the availability of gradients provide limited help
towards globally optimizing the multi-modal acquisition function. Global optimality guarantees are out of reach, since
branch-and-bound methods or DIRECT do not scale \citep{Jones1998}. Aside from local minima, another pitfall is the
presence of large plateaus where the acquisition function value is flat, while
local optima may be peaked. To decrease these two effects, \cite{rana2017high}
suggest to artificially inflate the gradient by taking a larger lengthscale
within an isotropic GP. Then the lengthscale is successively decreased while
tracking subsequent optima (warm started by previous values). Applying the random linear embedding of REMBO (Section \ref{sec:GP_AS})
at the level of acquisition function optimization, \citet{tran2019trading}
optimize over a finite set of subspaces.

Regardless of the quality of the fitted GP model, the expanding boundary volume as a function of the number of variables has a large impact on the acquisition's optimization (Equation \ref{eq:opt_acq}).
If based on uniform sampling, e.g., with multi-start gradient optimization or even evolutionary algorithms, the search will focus in these boundary areas. 
This is reinforced by the larger variance on the boundary that would locate the acquisition function's optimum there, with no hope of evaluating all the faces or vertices before focusing on the interior. 
Commenting on and discussing these effects, \citet{oh2018bock} proposed the use of cylindrical coordinates with BOCK to upweight the volume of the interior, imparting the prior knowledge
that the optimum is close to the center. Adding virtual derivative
observations on the boundary like in \cite{siivola2018correcting} is hardly
feasible in high dimension, as derivative information is only analytically enforceable at finitely many collocation points in GPs (and indeed computation scales cubically in the same). An infinite version could possibly be entertained via spectral methods, e.g., based on \cite{gauthier2012spectral}.
In such a case, the use of a trust-region (TR) can limit the
size of the search space drastically.
Trust region methods focus the optimization locally within the neighborhood of a current best solution (TR center), whose size increases if the new candidate point improves sufficiently over the TR center, or decreases instead, 
see e.g., \cite{larson2019derivative} for a review. 
Combined with BO, this has been shown to be quite
beneficial in high-dimension
\citep{regis2016trust,eriksson2019scalable,diouane2021trego,Zhou2021,daulton2021multiobjective}, perhaps at the
cost of a less global search (possibly compensated for with restarts or parallel
TR). To further avoid the attraction of the boundary of the TR,
(\ref{eq:opt_acq}) is only optimized over a discrete set in
\cite{Eriksson2018}, with some coordinates randomly kept at the TR center.

While the strategies above are mostly independent of the GP model (e.g., used
with isotropic or anisotropic product kernels), we next detail strategies for
GPs with structural assumptions. Note that these assumptions, never perfectly fulfilled
or estimated in practice, require the introduction of a noise
component to account for the introduced approximations, with a few exceptions. Hence this
has to be taken into account in the optimization of (\ref{eq:opt_acq}). 

\subsection{Additive case}

Additivity provides the opportunity to rely on the same decomposition as the
GP for solving (\ref{eq:opt_acq}), hence reducing the search to several
lower-dimensional searches, possibly in parallel. That is, rather than
focusing on the posterior of $f$, the idea is to look at those of the additive
components instead, the $g_i$ from models (\ref{mod:add}, \ref{mod:add_block},
\ref{mod:anova}): $\mathcal{N}(m_{n,i}(\x_I), s_{n,i}^2(\x_I))$ where
$s_{n,i}^2(\x_I) = k_i(\x_I, \x_I) - \veck_i(\x_i)^\top \K^{-1} \veck_i(\x_i)$
for a general index $I$. Then partial acquisition functions are defined on
each $g_i$ model, as a sum, e.g., as in \cite{Kandasamy2015}. More care is
needed for overlapping subsets, say relying on message passing
\citep{rolland2018high,hoang2018decentralized} but it remains more efficient
than optimizing in the original space.
The major advantage of this approach is that it allows for solving $d$ many 
one dimensional optimization problems instead of one $d$ dimensional problem, 
which \textit{ceteris paribus} is many times easier in the nonconvex case.
Nevertheless, the extent to which the
zero variance at unobserved locations with additive models affects the
globality of the search remains unknown. Adding a kernel a component that helps
mitigate incorrect assumptions as for variable selection or AS as below
may be interesting.
Additive GPs should be considered for problems where the high dimensionality
challenges the traditional acquisition function optimization and 
the budget is small enough that we cannot expect to learn a high-fidelity representation of the function.

\subsection{Variable selection}

If variables are completely removed beforehand in a preliminary stage, then a
fixed value may be used, returning to a low dimensional problem. Observations
with other values are generally discarded to keep the problem deterministic,
the alternative being to add some noise if the inert variables can vary.
Otherwise, after optimizing (\ref{eq:opt_acq}) on few variables, values for
the screened variables must be determined to evaluate $f$. Alternatives
include fixing these to a constant value, taking the values at the best design
sampled so far for these coordinates, a random sampling, or a combination of
these, see e.g., \cite{Li2017,spagnol2019global}. In \cite{Salem2018},
alternative lengthscales are estimated for the remaining variables by finding
the most different values still passing a likelihood ratio test. Then, values
for these variables are selected where the difference between the predictive
means are the most different between the two sets of hyperparameters, to
challenge the initial split.

When considering the acquisition function, the main advantage of variable selection approaches is to reduce the dimension
of the search space, which lessens the burden of the acquisition optimization.
On the other hand, the problem doesn't decompose, as in the additive case, and may be using suboptimal
values for the ``inactive" variables.

\subsection{Embedding case}
\label{sec:BO_AS}

With embeddings, additional difficulties arise with bounded domains. Let
$\mathbf{W} = [\A ~ \W_2]$ be a basis of $\R^d$. Splitting between active and
inactive (or less active) variables: $\forall \x \in \R^d$, $\x = \W
\W^\top \x = \A \A^\top \x + \W_2 \W_2^\top \x = \A \y + \W_2 \z$, $\y
\in \R^r$, $\z \in \R^{d- r}$. 
If $f$ has a true active subspace, the problem becomes 
\begin{equation}\text{find~} \y^* \in
\argmin \limits_{\y \in \mathcal{Y} \subseteq \R^r} f(p_\Xset(\A \y))
\label{eq:pb_REMBO}
\end{equation}
with $p_\Xset$ the convex projection onto $\mathcal{X}$; otherwise, the problem is: 
\begin{equation*}\text{find~} (\y^*, \z^*) \in \argmin \limits_{\y
\in \mathcal{Y} \subseteq \R^r, \z \in \mathcal{Z} \subseteq \R^{d-r}} f(\A \y + \W_2 \z) \text{~s.t.~} \A \y + \W_2 \z \in
\Xset.
\end{equation*} 
As discussed by \cite{Constantine2015}, perhaps the optimization over
$\z$ is 
secondary, which can be rewritten in the form $$\text{find~} \y^* \in
\argmin \limits_{\y \in \mathcal{Y} \subseteq \R^r} \min \limits_{\z \in \mathcal{Z} \subseteq \R^{d - r}} f(\A \y + \W_2
\z) \text{~s.t.~} \A \y \in \Xset$$ where less effort can be put on solving
the high-dimensional subproblem involving $\z$. This can be performed by
random sampling \citep{Constantine2015}. 
Alternatively, we can take a page from optimization
in the context of variable selection and fix the inactive directions, see e.g.,
\cite{cartis2020constrained}. Innocuous for unbounded domains\footnote{At
least in the formulation, unbounded domains are another topic in BO, see e.g.,
\cite{shahriari2016unbounded}.}, the constraint to belong to the domain for
evaluation is complex for a compact domain $\X$. 
The optimization problem is noisy whenever $\A$ changes over iterations as design points need to be projected on the embedding (that is, when a $\z$ component is present but ignored in the modeling of $f$).

Discussing solely the (centered) $[-1,1]^d$ hypercubic case, the intersection
with a linear embedding is a convex polytope, defined by $\mathcal{Y} = \left\{ \y \in \R^r
\text{~s.t.~} -1 \leq \A \y \leq 1 \right\}$. This is the domain advocated by
\cite{letham2020re} to solve (\ref{eq:pb_REMBO}), at the extra cost of
handling the linear constraints in optimizing $\alpha$. With random $\A$, a
simpler option is to optimize $\y$ within a hypercubic domain $\mathcal{Y} = \left[-l,
l\right]^r$, with $l$ chosen based on the probability of finding a
solution, as derived in \cite{Wang2013,Wang2016,qian2016random,cartis2020constrained}. These
results are simpler if the true active subspace is axis aligned, and depend on
the difference between $r$ and the true low dimension (if it exists). Nonetheless, these
choices for $\mathcal{Y}$ may not even contain $\y^*$. The smallest compact
set for this is a star-shaped polygon described in \cite{binois2020choice}.
The reason to focus on smaller $\Y$ is two-fold. First, most of the
difference corresponds to points on the boundary of $\Xset$, with increasingly
large volume compared to the intersection, bringing back the curse of
dimensionality. Second is that the corresponding observations are distorted
with the convex projection and are thus harder to model. These issues are
by-passed by \cite{nayebi2019framework,moriconi2020high}, at the cost of the
stricter diagonal or axis-aligned embedding assumptions.
Changing the domain, as is implicitly done by choosing $l$ relates to
another strategy not explored much in BO: to reduce the optimization space, as
discussed by \cite{Stork2020} or employed with TRs.

A weighted-PCA estimation of $\A$ is proposed in \cite{raponi2020high}, and they handle the domain issue with a penalty for infeasibility.
\cite{chen2020semi} use a semi-supervised version of sliced inverse regression (SIR) to find important input directions, using both labeled (evaluated) and unlabeled (unevaluated) designs. Namely, they
collect points which have high acquisition function values but were not chosen
for evaluation into an ``unlabeled dataset", which they incorporate into the
estimation of the embedding. The ``inverse"  in SIR indicates the goal of swapping the role of $\x$ and $y$, finding which values of the designs lead to given values of the outputs. Selecting a domain suffers the same difficulties as above, plus those related to updating the SIR model.

Most existing work fixes the embedding, or does not propagate the uncertainty when
estimated iteratively. Exploring the interplay between acquisition functions defined over AS
functions, e.g., in \cite{Garnett2013,wycoff2021gradient} and those defined for optimization
is a promising direction to enhance both of them. 

So when it comes to embeddings, they more often complicate than simplify the acquisition function search.
This can still be worth the extra effort when a sufficient evaluation budget exists to learn the embedding.

These problems of domain selection are amplified for non-linear dimension
reduction as illustrated by \cite{siivola2021good}. \cite{Li2016} faces
both additivity and embedding issues, restricting the extent of the projection
term to remain mostly in the domain. 
Next we provide pointers to benchmark problems to empirically assess the various methods detailed above.

\section{Synthetic test optimization problems}
\label{sec:tests}

Due to the multiplicity of possible methods and combinations, a thorough
empirical comparison is out of scope in this work. Besides, implementations are not always
available or even compatible. Inference techniques may differ as well, complicating comparisons. 
A less
obvious difficulty is the absence of standard benchmark functions for
high-dimensional surrogate-based optimization. We thus list here options we have come across.

First, some of the standard benchmark functions from
global optimization, say as in \cite{hansen2021coco}, can be scaled to high
dimension.
These synthetic test functions may be grouped into those which are separable, unimodal with moderate (resp.\ high) conditioning and finally multi-modal with  global structure. Only the last one is in the realm of BO, as stated by \cite{diouane2021trego}.
But even these are often pathological, for instance with exponentially
many local optima as $d$ increases, 
and consequently demand a greater budget than could be reasonably afforded to find the global solution, regardless of how quickly the simulator can be evaluated. Defining reasonable, sub-globally-optimal targets for these scalable problems under limited budgets could be an option.
Table \ref{tab:pbs_syn} of the Appendix enumerates these global optimization functions. Additional ones, created for BO, are tested by \cite{siivola2018correcting}.
The five physically-motivated analytical functions from \url{https://www.sfu.ca/~ssurjano/optimization.html} have been the subject of BO benchmarking too: Borehole ($8d$), OTL circuit ($6d$), piston ($7d$), robot arm ($8d$), and  wingweight ($10d$). The power circuit function ($13d$) \citep{Lee2019} makes one additional such example. 
Several functions designed for sensitivity analysis have been experimented on as well: the Sobol $g$-function (any $d$), the Ishigami function ($3d$), and the $8d$ flood model in \cite{Iooss2015}. Assessing the adequacy of these analytical examples with the structural model assumptions is interesting future work.

Another option is to add artificial, inactive variables to classical low
dimensional test functions, making the problem nominally high dimensional yet intrinsically low dimensional. While
this scales easily to billions of variables as in \cite{Wang2013}, it may not
be realistic. Repeated versions of these low dimensional toy functions are
also proposed, e.g., by \cite{oh2018bock}. This approach simply sums many low-dimensional versions of the same function together, allowing for more active variables.
Another option is to rotate the low-dimensional function via a linear transformation, with two snags: the initial hypercubic domain may not fill the entire high-dimensional one, or parts may be mapped outside, including known optima. 
Analytical functions can be extended easily to fill the gaps, but the difficulty of the problem becomes rotation dependent. 
A set of fixed rotations could be used for ease of comparison, or embedding matrices like the class proposed in \cite{nayebi2019framework} that do not have the same domain issues. 

GP realizations are an easy option to test the estimation difficulty of a given structure, and may also be used for comparisons between structures as in \cite{Ginsbourger2016}. 
The extent to which these results apply to real applications remains a topic of research.
To this end, a list of problems which may be reproduced using publicly available software is provided in Table \ref{tab:pbs}.
Some of these problems arise from hyperparameter tuning applications that are popular in the
machine learning community, and consist of tuning various neural networks
properties with BO. The choice of the neural networks structure seems to often be
hand tuned, the results are noisy and platform dependent, but the
code is sometimes available to reproduce results.
The most realistic examples comes from engineering and simulation, but the pipeline to reproduce high-fidelity simulation is rarely available in these cases.
Besides the ones in Table \ref{tab:pbs} for which this is the case, examples include the design of airfoils, wings or fans \citep{chen2020airfoil,gaudrie2020modeling,Palar2017,Viswanath2011,lukaczyk2014active,seshadri2019dimension}, automotive industry test cases \citep{Binois2015c}, alloy design \citep{Li2017,rana2017high}, biology \citep{ulmasov2016bayesian}, physics \citep{kirschner2019adaptive,mutny2018efficient}, and electronics \citep{Jones1998}.

\begin{table}[htpb]
    \caption{Examples of reproducible, realistic test problems.}
    \centering
		\begin{tabular}{ r  |c | c| l}
			Name(s) & Field/Domain & $d$  &  Example reference\\ 
			\hline
		    Bayesian neural net & Hyperparameter tuning & 5 & \cite{Falkner2018}\\
			Cartpole swing-up task & Hyperparameter tuning & 8 & \cite{Falkner2018}\\
			Covid-19 model & Epidemiology & 8 & \cite{wycoff2021sequential}\\
			Elliptic PDE  & Simulation & 100 & \cite{Tripathy2016}\\
			Groundwater remediation & Environment & 6 & \cite{Gramacy2020}\\
			Neural net & Hyperparameter tuning & 6 & \cite{Falkner2018}\\
			Robot pushing & Robotics & 14 & \cite{Wang2018}\\
			Rover trajectory & Robotics & 60 & \cite{Wang2018}\\
			Satellite drag & Astrodynamics & 7 & \cite{Gramacy2020}\\
			SDSS & Cosmology & 9 & \cite{Kandasamy2015}\\
			Shallow Net Weights & Hyperparameter tuning & 100-500 &
			\cite{oh2018bock}\\
			Three-link walk robot & Robotics & 25 & \cite{Zhang2019}\\
			VJ Cascade classifier & Hyperparameter tuning & 22 & \cite{Kandasamy2015}\\ %
		\end{tabular}

    \label{tab:pbs}
\end{table}

High dimensional data sets are another proper option to test high-dimensional GP regression, and can be extended to BO by allowing the BO algorithm to select not arbitrary points in the input space, but only those points corresponding to observed data. However, this discounts the effect of continuous optimization of the acquisition function, a fully nonneglible component of BO on real black-boxes. 
Still, some common data sets for regression tasks that have been used are listed in Table \ref{tab:pbs_dat}.
One option, following the example of \citep{jones2008large} presenting at MOPTA, is to fit an interpolant to a given dataset and then treat that interpolant as the black-box to be optimized.
We should remark that if the interpolant used to generate the ``black-box" is the same as that used for BO, the results may be unrealistically rosy.

\begin{table}[htpb]
    \caption{Example of data sets for regression used for GPs and BO.}
    \centering
		\begin{tabular}{ r | c| c | l}
			Name(s) & $d$ & $n$ & Example reference\\
			\hline 
			Airfoil & 5 & 1,503 & \cite{Francom2019}\\
			Amazon commerce reviews & 10,000 & 1,500 & \cite{Fukumizu2014}\\
			Autos &  24 & 159 &  \cite{delbridge2019randomly}\\
			Communities and crime & 96 & 1,993 & \cite{Garnett2013}\\
			Concrete & 8 & 1,030 &  \cite{Garnett2013}\\
			Crash & 9 & 20 &  \cite{Gramacy2020}\\
			CT slices & 318 & 3,071 & \cite{Garnett2013}\\
			Elevators & 17 & 8,752 &  \cite{gilboa2013scaling}\\
			Engine block and head joint sealing & 8 & 27 & \cite{Joseph2008}\\
			Heat exchanger & 4 & 64 & \cite{lin2020transformation}\\
			Housing & 13 & 506 & \cite{delbridge2019randomly}\\
			Isomap faces & 4,096 & 698 &  \cite{Guhaniyogi2016}\\
		    Kink40k & 8 & 10,000 & \cite{gilboa2013scaling}\\
		    LGBB & 3 & 3,167 &  \cite{Gramacy2020}\\
		    MARTHE & 20 & 300 &  \cite{Marrel2008}\\
		    Olivetti faces data set (modified) & 1,935 & 400 &  \cite{delbridge2019randomly}\\
		    PIC & 9 & 45,730 &  \cite{hoang2018decentralized}\\
			Pumadyn & 8 & 7,168 &  \cite{gilboa2013scaling}\\
			Pumadyn & 32 & 7,168 &  \cite{gilboa2013scaling}\\
			Sarcos robot & 21 & 44,484 & \cite{Winkel2021}\\
			Servo & 3 & 167 &  \cite{Duvenaud2011}\\
            Tecator & 100 & 215 &  \cite{plate1999accuracy}\\
            Temperature & 106 & 10,675 &  \cite{Garnett2013}\\
			Yacht & 6 & 308 &  \cite{Garnett2013}\\
		\end{tabular}
    \label{tab:pbs_dat}
\end{table}

From the set of examples listed above, we see that the definition of ``high dimensional" in BO varies extensively, from the low tens to several thousands. The optimization budgets share the same range of variation, but each analysis is conducted only for a fixed budget. There is thus room for a thorough study on which high-dimensional structure is the most adapted given the dimension and budget, or even better, enabling more complexity as more evaluations become available.

\section{Conclusion and perspectives}
\label{sec:conc}

Often rightly presented as one of the top challenges for GP-based BO, high-dimensionality has
generated many different ideas to address each aspect of the manifestation of
the curse of dimensionality. 
While the structural assumptions seem to cluster into one or more of variable selection, additive decomposition or linear
embeddings, no consensus is present on the best way to tackle the problem.
Even if this is largely a problem-dependent issue, the
inference methods used also restrict which structures are possible to assume, whether relying on a random structure or actually trying to infer it. Even besides these specifics, the particular strategy used for acquisition function optimization can overshadow
the rest. More systematic comparisons are necessary, as is the definition of
suitable benchmarks.

\subsection{Some general guidelines}

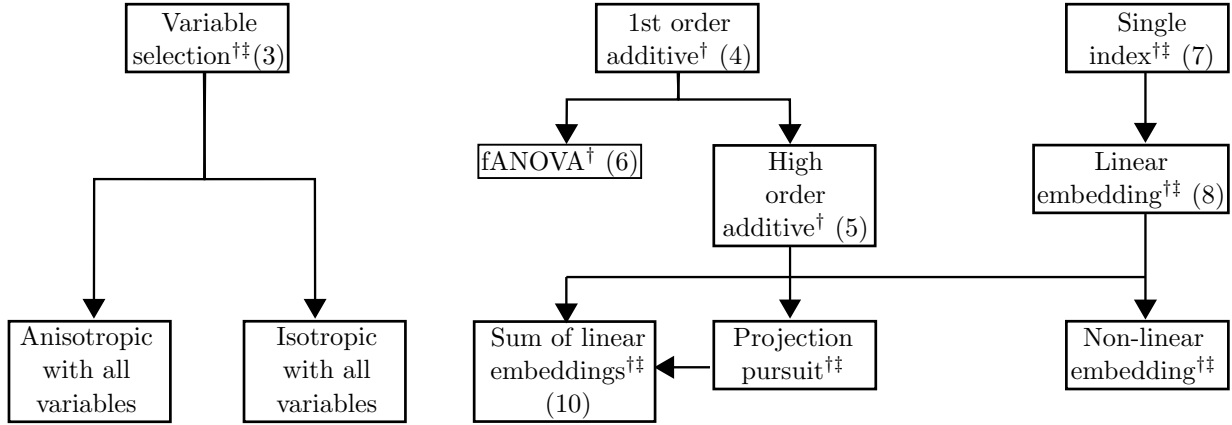
\begin{figure}[htpb]%
  \centering
  \def\svgwidth{\textwidth}
\begingroup%
  \makeatletter%
  \providecommand\color[2][]{%
    \errmessage{(Inkscape) Color is used for the text in Inkscape, but the package 'color.sty' is not loaded}%
    \renewcommand\color[2][]{}%
  }%
  \providecommand\transparent[1]{%
    \errmessage{(Inkscape) Transparency is used (non-zero) for the text in Inkscape, but the package 'transparent.sty' is not loaded}%
    \renewcommand\transparent[1]{}%
  }%
  \providecommand\rotatebox[2]{#2}%
  \newcommand*\fsize{\dimexpr\f@size pt\relax}%
  \newcommand*\lineheight[1]{\fontsize{\fsize}{#1\fsize}\selectfont}%
  \ifx\svgwidth\undefined%
    \setlength{\unitlength}{1473.85294276bp}%
    \ifx\svgscale\undefined%
      \relax%
    \else%
      \setlength{\unitlength}{\unitlength * \real{\svgscale}}%
    \fi%
  \else%
    \setlength{\unitlength}{\svgwidth}%
  \fi%
  \global\let\svgwidth\undefined%
  \global\let\svgscale\undefined%
  \makeatother%
  \begin{picture}(1,0.34418645)%
    \lineheight{1}%
    \setlength\tabcolsep{0pt}%
    \footnotesize
    \put(0,0){\includegraphics[width=\unitlength,page=1]{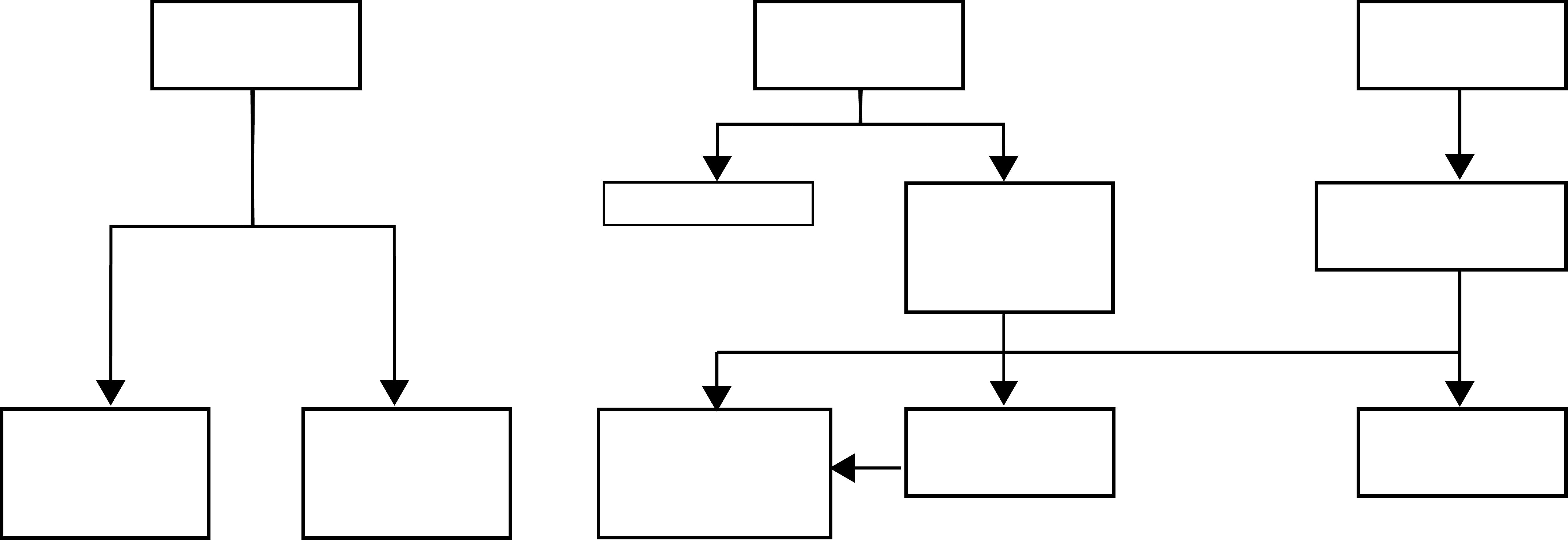}}%
    \put(0.585,0.2076043){\color[rgb]{0,0,0}\makebox(0,0)[lt]{\lineheight{1.25}\smash{\begin{tabular}[t]{c}High\\order\\additive$^{\dag}$ (\ref{mod:add_block})\\\end{tabular}}}}%
    \put(0.8425,0.20760437){\color[rgb]{0,0,0}\makebox(0,0)[lt]{\lineheight{1.25}\smash{\begin{tabular}[t]{c} Linear \\embedding$^{\dag\ddag}$~(\ref{mod:ridge})\\\end{tabular}}}}%
    \put(0.87,0.06335761){\color[rgb]{0,0,0}\makebox(0,0)[lt]{\lineheight{1.25}\smash{\begin{tabular}[t]{c}Non-linear\\embedding$^{\dag\ddag}$\\\end{tabular}}}}%
    \put(0.885,0.32102845){\color[rgb]{0,0,0}\makebox(0,0)[lt]{\lineheight{1.25}\smash{\begin{tabular}[t]{c}Single\\index$^{\dag\ddag}$ (\ref{mod:single})\end{tabular}}}}%
    \put(0.49,0.32102845){\color[rgb]{0,0,0}\makebox(0,0)[lt]{\lineheight{1.25}\smash{\begin{tabular}[t]{c}1st order\\additive$^{\dag}$ (\ref{mod:add})\end{tabular}}}}%
    \put(0.3875,0.20760432){\color[rgb]{0,0,0}\makebox(0,0)[lt]{\lineheight{1.25}\smash{\begin{tabular}[t]{c}fANOVA$^\dag$ (\ref{mod:anova})\\\end{tabular}}}}%
    \put(0.1025,0.32102845){\color[rgb]{0,0,0}\makebox(0,0)[lt]{\lineheight{1.25}\smash{\begin{tabular}[t]{c}Variable\\selection$^{\dag\ddag}$(\ref{eq:varsel}) \end{tabular}}}}%
    \put(0.01,0.06359682){\color[rgb]{0,0,0}\makebox(0,0)[lt]{\lineheight{1.25}\smash{\begin{tabular}[t]{c}Anisotropic\\with all\\variables\\\end{tabular}}}}%
    \put(0.22,0.06359682){\color[rgb]{0,0,0}\makebox(0,0)[lt]{\lineheight{1.25}\smash{\begin{tabular}[t]{c}Isotropic\\with all\\variables\\\end{tabular}}}}%
    \put(0.59365372,0.0633575){\color[rgb]{0,0,0}\makebox(0,0)[lt]{\lineheight{1.25}\smash{\begin{tabular}[t]{c}Projection\\pursuit$^{\dag\ddag}$\\\end{tabular}}}}%
    \put(0.395,0.06335761){\color[rgb]{0,0,0}\makebox(0,0)[lt]{\lineheight{1.25}\smash{\begin{tabular}[t]{c}Sum of linear\\embeddings$^{\dag\ddag}$\\ (\ref{mod:general})\\\end{tabular}}}}%
    \normalsize
  \end{picture}%
\endgroup%
\caption{A taxonomy of structural model assumptions, going from simpler (top) to more general ones (bottom). Arrows marks models that can be generalized, ($\dagger$) indicates that the estimation of a noise term is needed and ($\ddag$) cases for which a filling strategy is needed unless optimization is conducted only on a low dimensional manifold.}
\label{fig:taxo}
\end{figure}

As mentioned at the beginning of Section \ref{sec:highGP}, most structural model assumptions can be seen as instances of the more general model:
\begin{equation}
\text{model:~} f(\x) \approx \sum \limits_{i = 1}^\kappa g_i(\A_i \x)
\label{mod:general}
\end{equation}
with $\A_i \in \R^{r \times d}$, $\kappa \in \mathbb{N}^*$, which allows for sums of low dimensional effects, 
This includes at one end the additive Gaussian process of \cite{Durrande2010} which separates each variable individually into its own bin, at another end linear embeddings of the form $f(\x) = h(\mathbf{A} \x)$ which ignore certain input directions, and beyond lie the functional ANOVA kernels of \cite{Muehlenstaedt2012} or in combining several random low dimensional subspaces such as in \cite{delbridge2019randomly}. There are still promising research directions at the intersections of these methods, and we argue in particular that ``default'' kernel components, for instance a simple isotropic kernel, could be useful in complementing methods which search a particular part of the input space.
The various structural modeling options are presented in Figure \ref{fig:taxo}, highlighting the links between options.

Rather than starting directly from the most general model, whose direct inference would make a very steep initial step, we recommend to try the simpler models first. 
That is, begin with a standard GP to build an alternative one based on the ARD principle, or to compute sensitivity indices like Sobol indices.
Then build a fully additive GP model or a single index one as the basic models of the additive and linear embedding families.
From comparing the prediction given by these simple models, building more advanced models from the best performing one is possible: functional ANOVA or block additive from the additive one; linear embedding on top of single index. 
If none is performing well, projection pursuit or non-linear embeddings could be entertained, if $n$ is sufficiently large to allow a reliable inference of the corresponding model.

More precisely for the inference, depending on $d$, the budget $n$ and the application dependent complexity of $f$, inference by maximum likelihood (or, even better, a Bayesian alternative) is preferred. If too difficult or slow when $d$ reaches the hundreds, relying on a random structure to bypass the full-dimensional model inference remains possible. The simplest models scale well in computation with increasing $d$, while the more complex ones may benefit from recent advances brought by the increasing use of GPUs and automatic differentiation frameworks. This is especially true as more complex models require larger values of $n$ to show improvements.

Alternatively, if the budget of optimization is limited, a simple filling strategy for the ``inactive" variables combined with optimization on some restricted subspace could be preferable.
Finally the use of the trust region framework has shown promising results for optimization in high dimensional problems, and is amenable to parallelization.

\subsection{Future work}

There is potential for hybrids between models, inference methods and mathematical optimization techniques. For global optimization, a research avenue
pointed out by \cite{spagnol2019global} is to focus on variables that are important
to reach low function values, rather than evaluating their influence on the whole range of $f$. 
\cite{chen2020semi} go some way to focusing on low values, as do \cite{Guhaniyogi2016} on manifold-bound inputs. 
The ability to doubt the current structural assumption may also be
quite helpful in order to avoid negative feedback loop effects appearing in adaptive
design, see e.g., \cite{Gramacy2020}.

Whether the ability to keep a kernel defined on the whole initial domain is superior to defining one directly on a low dimensional embedding (which while appealing is associated with domain selection issues) remains to be seen.
If a kernel defined on an embedding is indeed preferred, keeping kernels sufficiently expressive (e.g., not restricted to be first order additive like in \cite{Durrande2013}) by adding an isotropic component may be improved by considering ortho-additivity \citep{Ginsbourger2016} instead. These orthogonality constraints are likely to help inference by reducing identifiability issues.

Additionally, the success of trust region BO approaches which restrict the search space to some subset of the input space motivates dimension reduction approaches which are themselves local. For instance \cite{WycoffThesis} propose to define an active subspace with respect to probability measures with non-uniform densities in order to emphasize certain regions of interest. Namely, the focus is restricted to a trust region, but other schemes could be considered. Of course, a dimension reduction which is \textit{locally} linear is in fact globally nonlinear, and such local linear models may prove to be a tractable approach to developing nonlinear dimension reduction for BO.

Furthermore, there have already been useful results born of the influence of mathematical programming on approaches to high dimensional Bayesian optimization. In particular, use of a rectangular trust region seems to be of benefit in \cite{eriksson2019scalable} and a Gaussian process analogue of the augmented Lagrangian method is developed in \cite{Gramacy2016,picheny2016bayesian}. Helping to develop new algorithms which combine benefits from low dimensional modeling with trust regions or other approaches successful in high dimensional mathematical programming broadly could help BO finally make the leap to practical, widespread usage on complex, high dimensional problems.


\appendix
\section{Synthetic benchmark functions}

Some global optimization benchmark functions are used for testing high-dimensional BO, of which a list is provided in Table \ref{tab:pbs_syn}.

\begin{table}[htpb]
    \caption{Synthetic problems for global optimization benchmarks tested on for high-dimensional BO. See e.g.,  \url{https://www.sfu.ca/~ssurjano/optimization.html} for further details, expressions and codes, or \cite{cartisdimensionality}. Cola is used by \cite{binois2020choice} while the Beach, Dream and Simba functions can be found in \cite{Winkel2021}.}
    \centering
		\begin{tabular}{ r | c|| r | c}
			Name(s) & $d$ & Name(s) & $d$ \\
			\hline
			Ackley & Any & Hartmann 3 & 3\\
			Beach & 6 & Hartmann 6 & 6 \\
			Beale & 2 & Levy & Any \\
			Branin & 2 & Michalewicz & Any\\
			Brent & 2 & Perm 4, 0.5 & 4 \\ 
			Bukin N.6 & 2 & Quadratic Form & Any \\
			Camel & 2 & Rastrigin & any \\
			Cola & 17 & Rosenbrock & Any\\
			Colville & 4  & Simba & 6 \\
			Dream & 6 & Shekel 5,7,10 & 4 \\
			Easom & 2 & Shubert & 2 \\
			Franke & 2 & Styblinski Tang & Any \\
			Friedman & 5  & Trid & any \\
			Giunta & 2 & Welsh & 20\\
			Goldstein-Price & 2 & Zakharov & any\\
			Griewank & Any & Zettl & 2 \\
		\end{tabular}
    \label{tab:pbs_syn}
\end{table}

\bibliographystyle{apalike}

\end{document}